\documentclass[a4paper]{amsart}
\usepackage[english]{babel}
\usepackage{amsfonts,amssymb}
\usepackage[bookmarks=false]{hyperref}
\newtheorem{theorem}{Theorem}

\newtheorem{lemma}[theorem]{Lemma}

\theoremstyle{definition}

\newtheorem{example}[theorem]{Example}
\newtheorem{remark}[theorem]{Remark}

\newtheorem*{scholium}{Scholium}
\newtheorem*{note}{Note}
\numberwithin{equation}{section}
\newcommand{\QQ}{\mathbf Q}
\newcommand{\RR}{\mathbf R}
\newcommand{\PP}{\mathbf P}
\newcommand{\GG}{\mathbf G}
\newcommand{\SL}{\mathbf{SL}}
\newcommand{\GL}{\mathbf{GL}}
\def\lra{\longrightarrow}
\def\hb{{\rm H}_{\rm b}}

\def\eI{{^{I}\!\mathrm{E}}}
\def\eII{{^{I\!I}\!\mathrm{E}}}
\def\dI{{^{I}\!d}}
\def\dII{{^{I\!I}\!d}}

\def\lft{L^\infty}
\def\am{$(\mathrm{A})$}
\def\hm{$(\mathrm{M}_\mathrm{I})$}
\def\hmm{$(\mathrm{M}_\mathrm{II})$}
\def\bu{\bullet}
\def\bsl{\backslash}
\begin{document}
\title[Vanishing up to the rank]{Vanishing up to the rank in bounded cohomology}
\author{Nicolas Monod}
\address{Universit\'e de Gen\`eve}
\email{nicolas.monod@math.unige.ch}
\thanks{Supported in part by Fonds National Suisse}
%\subjclass{}%
%\keywords{Bounded cohomology, Lie groups, algebraic groups}%
%\date{}%
%\dedicatory{}%
%\commby{}%
% ----------------------------------------------------------------
\begin{abstract}
We establish the vanishing for non-trivial unitary representations of the bounded cohomology of $\SL_d$ up to degree $d-1$. It holds more generally for uniformly bounded representations on superreflexive spaces. The same results are obtained for lattices. We also prove that the real bounded cohomology of any lattice is invariant in the same range.
\end{abstract}
\maketitle
\let\languagename\relax  % TO FIX A BUG IN RUNNING HEADERS AND BABEL
% ----------------------------------------------------------------
\section{Introduction}

By \emph{local field} we mean any non-discrete locally compact field, thus including the Archimedean as well as non-Archimedean cases of any characteristic. All unitary representations are assumed continuous and separable.

\bigskip

The cohomology of semisimple Lie groups and more generally algebraic groups over local fields is well-studied, see e.g.~\cite{Borel-Wallach} (we refer here to what is sometimes called \emph{continuous} group cohomology). It is particularly well understood for trivial coefficients~$\RR$ and for unitary representations. The cohomology of their \emph{lattices} is of importance since such lattices include the arithmetic groups. Thanks to Eckmann--Shapiro induction, the latter study can be reduced to the former, though only at the cost of considering more general coefficients in the non-uniform case; for this reason, the cohomology of arithmetic groups still presents difficulties.

\medskip

One important result for the cohomology of semisimple groups is \textbf{(i)}~the \emph{vanishing below the rank}, due to Borel--Wallach~\cite{Borel-Wallach}, Zuckerman~\cite{Zuckerman}, Casselman~\cite{Casselman} (see~V.3.3 and XI.3.9 in~\cite{Borel-Wallach}). For arithmetic groups, an interesting issue is \textbf{(ii)}~in what range their cohomology is \emph{invariant}, namely arises by restriction from the cohomology of the ambient semisimple group. The answer turns out to depend on the type of lattice, reflecting difficulties in applying induction for non-uniform lattices; for results depending notably on the $\QQ$-rank of arithmetic groups, see Borel~\cite[Theorem~7.5]{Borel74} (see also~\cite{Matsushima},\cite{Garland},\cite{Borel-Serre73}).

\medskip

This note investigates the analogue of the these two issues~(i),~(ii) for \emph{bounded cohomology~$\hb^\bu$}; we restrict ourselves to the groups $\SL_d$.

For a sampling of the uses of bounded cohomology, background and definitions, see~\cite{Gromov},\cite{Burger-Monod3},\cite{MonodICM}. Vanishing in degree two for groups of rank at least two was established in~\cite{Burger-Monod1} (see also~\cite{Burger-Monod3},\cite{Monod-Shalom1}; partial results in degree three for trivial coefficients are obtained in~\cite{Burger-MonodERN},\cite{MonodJAMS}). Above degree two, not much is known about the bounded cohomology of semisimple groups and their lattices (nor for any other group, for that matter).

\begin{theorem}
\label{thm:sl}%
Let $k$ be a local field,  $G=\SL_d(k)$ and $V$ a unitary representation of $G$ not containing the trivial one. Then
$$\hb^n(G, V)=0 \kern1cm\text{for all \ $0\leq n \leq d-1$.}$$
The same holds more generally for uniformly bounded continuous representations on superreflexive separable Banach spaces without invariant vectors.
\end{theorem}

As with ordinary cohomology, induction sometimes allows to deduce results about lattices from results on the ambient group. However, it takes us a priori to non-separable and non-continuous modules when dealing with bounded cohomology. It is only in degrees up to two that it has been possible to use unitary induction (\cite[Corollary~11]{Burger-Monod3}, \cite[11.1.5]{Monod}). Nevertheless, we prove:

\begin{theorem}
\label{thm:sl:lat}%
Let $k$ be a local field,  $G=\SL_d(k)$, $\Gamma<G$ a lattice and $W$ a unitary representation of $\Gamma$ not containing the trivial one. Then
$$\hb^n(\Gamma, W)=0 \kern1cm\text{for all \ $0\leq n \leq d-1$.}$$
The same holds more generally for uniformly bounded continuous representations on superreflexive separable Banach spaces without invariant vectors.
\end{theorem}

Whilst the bounded cohomology with \emph{trivial} coefficients remains mysterious for both $G$ and its lattices, we can still determine the relation between the two; observe that here the answer does not depend upon the type of the lattice:

\begin{theorem}
\label{thm:lattice}%
Let $k$ be a local field,  $G=\SL_d(k)$ and $\Gamma<G$ a lattice. Then the restriction
$$\hb^n(G, \RR) \lra \hb^n(\Gamma, \RR)$$
is an isomorphism for all \ $0\leq n \leq d-1$.
\end{theorem}

The three theorems above depend on one general statement, Theorem~\ref{thm:vr} below. In order to apply the latter, we need to understand the intersection of several conjugates of a suitable maximal parabolic subgroup $Q<G$. To decide whether vanishing results as above can be deduced from Theorem~\ref{thm:vr} for other classical simple groups, one would need to control the combinatorics of such intersections in general.

\begin{scholium}
One (too rare) instance in which bounded cohomology behaves better than ordinary cohomology is that the restriction to lattices is always injective~\cite[8.6.2]{Monod}. Therefore, one could deduce Theorem~\ref{thm:sl} from Theorem~\ref{thm:sl:lat} using Moore's theorem (recalled below). We found it more natural to begin by establishing the simpler Theorem~\ref{thm:sl}.
\end{scholium}

We conclude this introduction by observing that the bound on $n$ is sharp at least when $d=2$ for Theorems~\ref{thm:sl} and~\ref{thm:sl:lat}. Indeed, for $G=\SL_2(k)$ and any lattice $\Gamma$ the spaces $\hb^2(G, L^2(G))$ and $\hb^2(\Gamma, \ell^2(\Gamma))$ are non-zero (see e.g.~\cite{Monod-Shalom1}).

\section{General setting}

Let $G$ be a locally compact second countable group. A Banach space $V$ with an isometric linear representation is called a \emph{coefficient $G$-module} if it is the dual of a separable Banach space with continuous $G$-representation~\cite[1.2.1]{Monod}. Examples include unitary representations and $\lft$ spaces, the latter being in general neither continuous nor separable.

\begin{remark}
Theorems~\ref{thm:sl} and~\ref{thm:sl:lat} are stated for \emph{uniformly bounded} representations on superreflexive spaces. (We recall that one characterisation of superreflexivity is the existence of an equivalent uniformly convex norm~\cite[Theorem~A.6]{Benyamini-Lindenstrauss}.) Whilst one could perfectly well work with bounded cohomology in uniformly bounded representations, we will not need to. Indeed, one can always replace the norm with an equivalent invariant norm; since superreflexivity is a topological property, it is preserved under this operation. Thus we shall from now on assume all representations isometric.
\end{remark}

Let $Q<G$ a closed subgroup, $N\triangleleft Q$ an amenable closed normal subgroup and $r$ a positive integer. Endow $B=G/Q$ with its canonical invariant measure class (Theorem~23.8.1 in~\cite{Simonnet}). Consider the following conditions on intersections of generic conjugates of $Q$:

\medskip

\begin{itemize}
\item[\hm] For a.e. point in $B^{r-1}$, its stabiliser in $N$ has no non-zero invariant vectors in $V$.

\item[\hmm] For a.e. point in $B^{r+1}$, its stabiliser in $G$ has no non-zero invariant vectors in $V$.

\item[\am] The $G$-action on $B^{r+1}$ is amenable in Zimmer's sense~\cite{Zimmer78b},\cite[4.3.1]{Zimmer84}.
\end{itemize}

\smallskip

\noindent
(M)~stands for Moore. The motivating examples for this axiomatics are provided in the setting of Theorems~\ref{thm:sl},~\ref{thm:sl:lat} and~\ref{thm:lattice} (with $r=d-1$), see Example~\ref{exo:sl} below. Here is the main technical result of this note:

\begin{theorem}
\label{thm:vr}%
(i) If \hm\ holds, then $\hb^n(G, V)=0$ for all $0\leq n \leq r-1$.

\smallskip
\noindent
(ii) If in addition \hmm\ holds, then $\hb^n(G, V)=0$ for all $0\leq n \leq r$.

\smallskip
\noindent
(iii) If in addition \am\ holds, then the complex
$$0 \lra \lft(B, V)^G \lra \lft(B^2, V)^G  \lra \lft(B^2, V)^G \lra\cdots$$
realises $\hb^n(G, V)$ for all $n\geq 0$.
\end{theorem}

In the above complex and everywhere below, the maps are the usual \emph{homogeneous differentials} defined by $d=\sum (-1)^k d_k$, where $d_k$ omits the $k$-th variable.

\begin{remark}
\label{rem:coboundaries}%
Notice that $\hb^{r+1}(G, V)$ is the space of \emph{cocycles} on $B^{r+2}$ in case~(iii), since $\lft(B^{r+1}, V)^G$ vanishes by \hmm\ and thus there are no coboundaries. This implies notably that $\hb^{r+1}(G, V)$ is Hausdorff.
\end{remark}

The following classical result is essentially due to Moore.

\begin{theorem}
Let $k$ be a local field, $\GG$ a connected simply connected $k$-simple $k$-group and $V$ a unitary $\GG(k)$-representation not containing the trivial one. Then unbounded subgroups have no non-zero invariant vectors. The same holds more generally for uniformly bounded continuous representations on superreflexive Banach spaces without invariant vectors.
\end{theorem}

\begin{proof}
The most classical form is due to Moore~\cite[Theorem~1]{Moore66}. For general fields, see Proposition~5.5 in~\cite{Howe-Moore}. Howe--Moore further established the stronger statement that matrix coefficients vanish at infinity (see e.g.~\cite{Margulis},\cite{Zimmer84}). It was noticed by Shalom that the latter still holds for continuous isometric representations on uniformly convex uniformly smooth Banach spaces. For a proof of this fact, see the Appendix of~\cite{Bader-Furman-Gelander-Monod}. It remains only to justify that all uniformly bounded representations on superreflexive Banach spaces admit an equivalent uniformly convex and uniformly smooth norm which is invariant; this is proved in~\cite[Proposition~2.3]{Bader-Furman-Gelander-Monod}.
\end{proof}

\begin{example}
\label{exo:sl}%
Let $G=\SL_{r+1}(k)$, $B$ the projective space $\PP^r(k)$ and $Q\cong k^r\rtimes \GL_r(k)$ a corresponding maximal parabolic, $N\cong k^r$ its unipotent radical. One verifies directly that any conjugate of $Q$ in $G$ intersects $N$ along a subspace of codimension at most one. It follows that the stabiliser in $N$ of every (not a.e) point in $B^{r-1}$ is unbounded in $G$. One also checks that the intersection of any $r+1$ conjugates of $Q$ contains a conjugate of the diagonal subgroup of $G$, hence is also unbounded. Thus \hm\ and \hmm\ follow from Moore's theorem if $V$ is a uniformly bounded continuous representation on a superreflexive Banach space without invariant vectors.

As for \am, it is enough by~\cite[Theorem~A]{Adams-Elliott-Giordano} to verify that generic stabilisers for $B^{r+1}$ are amenable (in view of the smoothness of the $G$-action). And indeed, away from a finite union of subvarieties of positive codimension, these stabiliser are conjugated to the diagonal subgroup, which is amenable.
\end{example}

The above example shows that Theorem~\ref{thm:sl} follows from Theorem~\ref{thm:vr} with $r=d-1$. Notice that for the vanishing of Theorem~\ref{thm:sl}, condition \am\ is not needed; however, it provides additional information via Remark~\ref{rem:coboundaries}. In order to address Theorem~\ref{thm:lattice}, we need a variation on that example:

\begin{example}
\label{exo:lattice}%
In Example ~\ref{exo:sl}, the $G$-modules $V$ were separable and continuous. However, the same reasoning applies to a special class of further coefficient modules, namely $V=\lft_0(X)$ for any ergodic measure-preserving standard probability $G$-space ($\lft_0$ denotes the subspace of functions of integral zero). Indeed, Moore's theorem applied to $L^2(X)$ shows that any unbounded subgroup of $G$ still acts ergodically on $X$ and thus has no non-zero invariant vectors in $V$.
\end{example}

\begin{proof}[Proof of Theorem~\ref{thm:lattice}]
The restriction map is always injective~\cite[8.6.2]{Monod}. Moreover, it fits in an exact sequence~\cite[10.1.7]{Monod} in such a way that we just need to prove the vanishing of $\hb^n(G,-)$ for $\lft(G/\Gamma)/\RR$. The latter is isomorphic to $\lft_0(G/\Gamma)$, so that we can apply Theorem~\ref{thm:vr} as in Example~\ref{exo:lattice}.
\end{proof}

Finally, a combination of both approaches:

\begin{proof}[Proof of Theorem~\ref{thm:sl:lat}]
Let $V$ be the \emph{induced} $G$-module associated to $W$~\cite[10.1.1]{Monod}. That is, $V$ is the coefficient $G$-module $V=\lft(G,W)^\Gamma$ of (left) $\Gamma$-equivariant $\lft$-maps $f:G\to W$ endowed with the $G$-action by right translations on $G$. We claim that $V$ satisfies \hm\ and \hmm. Indeed, $V$ embeds into the space $U=L^{[2]}(G,W)^\Gamma$ of $\Gamma$-equivariant maps such that the function $\|f\|_W: \Gamma\bsl G\to \RR$ is $L^2$; that is, as a Banach space, $U\cong L^2(\Gamma\bsl G, W)$. In the unitary case, $U$ is still a Hilbert space with unitary representation. If $W$ is merely superreflexive, it is a result of Figiel--Pisier that $U$ is so too, see~\cite[II~1.e.9~(i)]{Lindenstrauss-Tzafriri}. Therefore, applying Moore's theorem to $U$, we deduce that the only elements of $V$ fixed by an unbounded subgroup are $G$-fixed. Since a $G$-fixed function class has its essential range in $W^\Gamma$ by the definition of $V$, it vanishes~--- proving the claim. Now Theorem~\ref{thm:vr} applies and $\hb^n(G,V)=0$. We conclude by the induction isomorphism~\cite[10.1.3]{Monod} which identifies the latter to $\hb^n(\Gamma, W)$.
\end{proof}

\section{Proof of Theorem~\ref{thm:vr}}

Define a first quadrant double complex by $L^{p,q} = \lft(G^{p+1} \times B^{q+1}, V)^G$. We will freely use the identifications
$$L^{p,q}\ \cong\ \lft\big(G^{p+1}, \lft(B^{q+1}, V)\big)^G\ \cong\  \lft\big(B^{q+1}, \lft(G^{p+1}, V)\big)^G$$
which follow e.g. from the Dunford--Pettis theorem, see~\cite[2.3.3]{Monod}. Define $\dI :L^{p,q}\to L^{p+1, q}$ by the homogeneous differential associated to $G^{p+1}$ and $\dII :L^{p,q}\to L^{p, q+1}$ by the homogeneous differential on $B^{q+1}$ affected with the sign $(-1)^{(p+1)}$. To such a complex are associated two spectral sequences $\eI,\, \eII$ defined respectively by
$$\eI_1^{p,q} = {\rm H}^{p,q}(L^{p,\bu}, \dII),\kern.7cm \eII_1^{p,q} = {\rm H}^{q,p}(L^{\bu,p}, \dI)$$
both abutting to the cohomology of the associated total complex.

\begin{note}
For background on spectral sequences we refer to~\cite[III\S14]{Bott-Tu},\cite[III.7]{Gelfand-Manin}. A similar bicomplex was used in~\cite{MonodJAMS}; we point out at this occasion that the induction step in that paper (over the rank $n$) needs to be changed: It holds in degrees $q\leq 2+n/2$ and proceeds separately over even and odd integers $n$.
\end{note}

\begin{lemma}
\label{lem:eI}%
The first spectral sequence converges to $\hb^\bu(G, V)$.\\
More precisely, $\eI_t^{p,q}=0$ $(\forall\,p\geq0, q\geq1, t\geq1)$ and $\eI_t^{p,0} \cong \hb^p(G, V)$ $(\forall\,p\geq0, t\geq2)$.
\end{lemma}

\begin{proof}
The cohomology of the complex
$$0\lra \lft(B,V) \lra \lft(B^2,V) \lra \lft(B^3,V) \lra \cdots \leqno{(*)}$$
is concentrated in degree zero, where it is $V$. Indeed, the augmented complex 
$$0\lra V \lra \lft(B,V) \lra \lft(B^2,V) \lra \lft(B^3,V) \lra \cdots$$
is acyclic, a contracting homotopy being provided by (Gelfand--Dunford) integration over $B$ of the first variable (as in~\cite[7.5.5]{Monod}, except here we do not need to worry about the continuous submodules). The functor $\lft(G^{p+1}, -)^G$ is exact for all $p\geq0$ with respect to \emph{adjoint} short exact sequences of coefficient $G$-modules~\cite[8.2.5]{Monod}. Applying it to~$(*)$ (which is adjoint~\cite[7.5.4]{Monod}), it follows that the cohomology of the complex $L^{p,\bu}$, namely $\eI_1^{p,q}$, is concentrated in degree zero, where it is $\lft(G^{p+1},V)^G$. This shows at once the vanishing for $q\geq 1$ and the isomorphism  $\eI_2^{p,0}\cong \hb^p(G, V)$ since the latter is realised by the complex $\lft(G^{p+1},V)^G$. (This can be taken as a definition of bounded cohomology; compare~\cite[7.5.1]{Monod}.)
\end{proof}

\begin{lemma}
\label{lem:eII_1}%
There is a canonical identification $\eII_1^{p,q}\ \cong\ \hb^q\big(G, \lft(B^{p+1}, V)\big)$ for all $p,q\geq0$.
\end{lemma}

\begin{proof}
The complex $\lft\big(G^{\bu+1}, \lft(B^{p+1}, V)\big)^G$ indeed computes $\hb^q\big(G, \lft(B^{p+1}, V)\big)$.
\end{proof}

\begin{lemma}
\label{lem:van_low}%
If  \hm\ holds, then $\eII_1^{p,q}=0$ for all $q$ and all $p\leq r-1$.
\end{lemma}

\begin{proof}
Consider the identification $\lft(B^{p+1}, V) \cong \lft(G/Q, \lft(B^p, V))$ (\cite[2.3.3]{Monod}). The induction isomorphism~\cite[10.1.3]{Monod}, the characterisation of the induction of restricted modules~\cite[10.1.2(v)]{Monod} and Lemma~\ref{lem:eII_1} imply that $\eII_1^{p,q}$ is isomorphic to $\hb^q\big(Q, \lft(B^p, V)\big)$. In the latter, we may replace the coefficients by $\lft(B^p, V)^N$ because the group $N$ is amenable~\cite[8.5.3]{Monod}. But condition \hm\ implies that the stabiliser in $N$ of a.e. point in $B^p$ has no non-zero invariant vectors in $V$ whenever $p\leq r-1$. Thus $\lft(B^p, V)^N$ vanishes and hence $\eII_1^{p,q}=0$.
\end{proof}

At this point, part~(i) of the theorem is established, since the cohomology of the total complex in degree $n$ involves only terms $\eII_\bu^{p,q}$ with $p+q=n$, which all vanish when $n\leq r-1$. For part~(ii), it suffices now to show that $\eII_1^{r,0}$ vanishes. By Lemma~\ref{lem:eII_1}, this amounts to $\lft(B^{r+1}, V)^G=0$, which follows from \hmm.

\begin{lemma}
If \am\ holds, then $\eII_1^{p,q}=0$ for all $q>0$ and all $p\geq r$.
\end{lemma}

\begin{proof}
Condition \am\ implies amenability for $B^{p+1}$ whenever $p\geq r$ (\cite[4.3.4]{Zimmer84}). Therefore, the $G$-module $\lft(B^{p+1}, V)$ is relatively injective~\cite[5.7.1]{Monod}. This implies that the right hand side in Lemma~\ref{lem:eII_1} vanishes~\cite[7.4.1]{Monod}.
\end{proof}

Thus in the setting of part~(iii) it follows that $\eII_t^{p,q}$ stabilises already at the tableau $t=2$ and that the cohomology of the total complex is just $\eII_2^{\bu,0}$, which is precisely computed on the complex
$$\eII_1^{\bu,0}:\ 0 \lra \lft(B, V)^G \lra \lft(B^2, V)^G  \lra \lft(B^2, V)^G \lra\cdots$$
as was to be shown.

% ----------------------------------------------------------------

\end{document}